\documentclass{elsart}
\usepackage{amsmath,amscd,amsfonts,amssymb}

\newtheorem{Theorem}{Theorem}[section]
\newtheorem{Lemma}[Theorem]{Lemma}
\newtheorem{Corollary}[Theorem]{Corollary}
\newtheorem{Proposition}[Theorem]{Proposition}
\newtheorem{Remark}[Theorem]{Remark}
\newtheorem{Example}[Theorem]{Example}

\def\deg{\operatorname{deg}}
\def\reg{\operatorname{reg}}
\def\gin{\operatorname{in}}
\def\codim{\operatorname{codim}}

\def\sk{\smallskip}
\def\Nset{{\mathbb N}}
\def\Zset{{\mathbb Z}}
\def\Rset{{\mathbb R}}

\def\mfr{{\mathfrak m}}
\def\abf{{\mathbf  a}}
\def\bbf{{\mathbf  b}}
\def\mbf{{\mathbf  m}}
\def\nbf{{\mathbf  n}}
\def\pbf{{\mathbf  p}}
\def\qbf{{\mathbf  q}}
\def\ebf{{\mathbf  e}}
\def\xbf{{\mathbf  x}}
\def\ybf{{\mathbf  y}}
\def\tbf{{\mathbf  t}}

\def\Pcal{{\mathcal P}}
\def\Acal{{\mathcal A}}

\begin{document}

\begin{frontmatter}

\title{ Gr\"obner bases of simplicial toric ideals}
\thanks{The second  author was supported by NAFOSTED 
(Vietnam) and Max-Planck Institute for Mathematics in the Sciences (Germany). }
\author{Michael Hellus}
\address{Universit\"at Leipzig, Fakult\"at f\"ur Mathematik und Informatik,
Augustusplatz 10/11, D-04109 Leipzig, Germany }
\ead{Michael.Hellus@math.uni-leipzig.de}
\author{L\^e Tu\^an Hoa}
\address{Institute of Mathematics Hanoi, 18 Hoang Quoc Viet Road, 10307 Hanoi, 
Vietnam}
\ead{lthoa@math.ac.vn}
\author{J\"urgen St\"uckrad}
\address{Universit\"at Leipzig, Fakult\"at f\"ur Mathematik und Informatik,
Augustusplatz 10/11, D-04109 Leipzig, Germany }
\ead{stueckrad@math.uni-leipzig.de}
\begin{abstract}  Bounds for the  maximum degree of  a minimal  Gr\"obner basis  of simplicial
toric ideals with respect to the reverse lexicographic order are given. These bounds are close to the bound stated in Eisenbud-Goto's
Conjecture on the Castelnuovo-Mumford regularity.
\end{abstract}

 \maketitle

\end{frontmatter}

\section*{Introduction}

Let $I$ be a homogeneous ideal of a polynomial ring $R$. The coarsest measure of the 
complexity of a Gr\"obner basis (w.r.t. to a term order $\leq $) of an ideal $I$ is its maximum 
degree, which is the highest degree of a generator of the initial ideal $\gin_\leq(I)$. However, 
this quantity is not easy to be handled with. One way to study it is to use a better-behaved 
invariant, the Castelnuovo-Mumford regularity $\reg (I)$ of $I$. This invariant can be defined 
as the maximum over all $i$ of the degree minus $i$ of any minimal $i$-th syzygy of $I$, 
treating generators as 0-th syzygies. In the generic coordinates and with respect to the reverse 
lexicographic order, the maximum degree in  a minimal Gr\"obner basis of $I$ is bounded by $\reg (I)$  
(see \cite[Corollary 2.5]{BS}). Unfortunately, this is not true for arbitrary coordinates (see, e.g., the example after \cite[Lemma 14]{HHy}).  On the
other hand, a famous conjecture by Eisenbud and Goto states that $\reg(I)\leq \deg (R/I) 
-\codim(R/I) + 1$, provided $I$ is a prime ideal containing no linear form (see \cite{EG}). Here 
$\deg(R/I)$ and $\codim(R/I)$ denote the multiplicity and the codimension of $R/I$, 
respectively. Thus, in the generic coordinates the {\it Eisenbud-Goto bound} $\deg (R/I) -\codim(R/I) 
+ 1$ is an expected bound for the maximum degree in a minimal  Gr\"obner basis of $I$ w.r.t. the reverse lexicographic order of a prime ideal containing no linear form. We may hope that this expectation still holds for some other coordinates. 

In this paper we are interested in estimating the degree-complexity of  Gr\"obner bases of  
simplicial toric ideals. Toric ideals are nice, particularly because they are prime ideals  and in the natural coordinates they are
generated by binomials.  In order to find a minimal Gr\"obner basis of such an  ideal,  it is therefore natural to try to keep the original coordinates, so that elements  of  such a Gr\"obner basis can be taken as binomials - which are cheap to compute and to restore.  On the  
other hand in \cite{HS} the last two authors have shown that for a large class of simplicial toric ideals $I$, the Castelnuovo-Mumford regularity $\reg (I)$ is bounded by the {\it Eisenbud-Goto bound} $\deg (R/I) -\codim(R/I)  
+ 1$.  From these phenomena  we  
believe that following conjecture holds:

 {\bf Conjecture}: {\it Assume that $I$ is the toric  ideal associated with a homogeneous simplicial 
 affine semigroup $S$ over an arbitrary  field $K$. The maximum degree in a minimal  Gr\"obner basis of $I$ in the natural coordinates and w.r.t. the reverse lexicographic order is bounded above by  $\deg K[S] -\codim K[S] + 1$.}

Note that this is not true for an arbitrary term order (see Example \ref{A1A3}). For the rest of the paper, if not otherwise stated,  we consider only the natural coordinates and the reverse lexicographic order.  Although we are still not able to solve the  
above problem, we can establish the upper bound $2(\deg K[S] -\codim K[S])$.  In order to do that we first establish an upper bound  in terms of the reduction number $r(S)$ of $K[S]$. Then, combining with a bound of \cite{HS} on $r(S)$, we get the main  result, see Theorem \ref{A1}. We also provide another bound  
 in terms of the codimension $c= \codim K[S]$ and the total degree $\alpha $ of monomials  
defining $S$ (Theorem \ref{A4}). In a lot of  examples  bounds in Theorems \ref{A1} and \ref{A4} are even much  
smaller then the Eisenbud-Goto bound. 

 In Section \ref{B} we  solve the above conjecture for certain classes of simplicial toric  
ideals. Ideals of first type come from a simple observation that the maximum degrees in their minimal  
Gr\"obner bases are bounded by the Castelnuovo-Mumford regularity if the corresponding rings   
$K[S]$ are  generalized Cohen-Macaulay rings. Ideals of second type are raised  by certain properties of  the parameter set $\Acal$ (see Propositions \ref{B2} and \ref{B3}).  In this situation, by using Theorem  
\ref{A4} we can restrict ourselves to  few exceptional cases when the codimension is very big.  
Then the main technique is  to refine bounds on the reduction number or to calculate its exact  
value, so that one can apply Theorem \ref{A1}. In particular,  we show that  the  conjecture holds for all simplicial toric  
ideals in \cite{HS}, for which the Eisenbud-Goto conjecture is known to be true.
\sk

{\it Notation}: In this paper we use bold letters to denote a vector, while their coordinates are  
written in the normal style. Thus $a_i,\ e_{1i}$ are the $i$-th coordinates of vectors $\abf,\  
\ebf_1$, respectively; $\xbf^{\mbf} = x_1^{m_1}\cdots x_c^{m_c}$, $\ybf^{\nbf} =  
y_1^{n_1}\cdots y_d^{n_d}$ and $\tbf^{\nbf} = t_1^{n_1}\cdots t_d^{n_d}$. The  ordering  
of variables is always assumed to be $x_1 > \cdots > x_c >y_1>\cdots > y_d$. We always write a binomial in such a way that its first term  is bigger  than the second one.

\section{Bounds}\label{A}

Let $S\subseteq \Nset^d$ be a homogeneous, simplicial affine semigroup generated by 
  a set of elements of the following type:
$$\Acal = \{\ebf_1,...,\ebf_d, \abf_1,...,\abf_c\} \subseteq M_{\alpha ,d} =\{(x_1,...,x_d)\in   
\Nset^d|\ x_1 +\cdots + x_d =\alpha \},$$
where $c\geq 2, \alpha \geq 2$ are natural numbers and $ \ebf_1 = (\alpha, 0,...,0), ..., \ebf_d  
=(0,...,0, \alpha)$.
Moreover, if $\abf_i = (a_{i1},...,a_{id})$, we can assume that the integers $a_{ij}$, where $  
i=1,...,c ,\ j=1,...,d$, are
relatively prime.  Note that $\dim K[S] = d$ and $\codim K[S]=c$.  Let $I_\Acal$ be   the  
kernel of the homomorphism
$$\begin{array}{r}K[\xbf,\ybf] := K[x_1,...,x_c,y_1,...,y_d] \rightarrow K[S]\equiv   
K[t_1^\alpha ,...,t_d^\alpha, \tbf^{\abf_1},...,\tbf^{\abf_c}] \subseteq K[\tbf];\\ x_i \mapsto  
\tbf^{\abf_i};\ y_j \mapsto t_j^\alpha,\ i=1,...,c;\ j=1,...,d .
\end{array}$$
We call $I_\Acal$ a {\it simplicial toric ideal defined by $\Acal$} (or $S$). We will consider  
the standard grading on $K[\xbf,\ybf]$ and $K[S]$, i.e. $\deg(x_i) = \deg(y_j) =1$ and if  
$\bbf\in S$, then $\deg(\bbf) = (b_1+\cdots + b_d)/\alpha$.

Note that  $I_\Acal$ always has a minimal Gr\"obner basis consisting of   
binomials (see, e.g.,  
\cite[Chapter 1]{St1}). We are interested in bounding  its  maximum degree.

Let $A=A_0\oplus A_1 \oplus \cdots$, where $A_0 = K$, be a standard graded
$K$-algebra of dimension $d$. A minimal reduction of $A$ is a graded ideal $I$
generated by $d$ linear forms such that $[IA]_n = A_n$ for $n\gg 0$. The least
integer $n$ such that $[IA]_{n+1} = A_{n+1}$ is called the reduction number of
$A$ w.r.t. $I$ and will be denoted by $r_I(A)$. Note that $(t_1^\alpha,...,
t_d^\alpha)$ is a minimal reduction of $K[S]$. We denote by $r(S)$  the reduction number
 of $K[S]$ w.r.t. this minimal reduction. Then $r(S)$ is the least positive
integer $r$ such that  $(r+1)\Acal = \{e_1,...,e_d\}+r\Acal$, where for two subsets $B$
and $C$ of $\Zset^d$ we denote by $B\pm C$  the 
set of all elements of the form $b\pm c,\ b\in B, \ c\in C$, and $nB = B+\cdots + B$ ($n$ times).   
This reduction number was used in \cite{HS} to bound the Castelnuovo-Mumford regularity of  
$K[S]$. 

\begin{Theorem} \label{A1} The maximum degree in a minimal  Gr\"obner basis of $ I_\Acal$ is bounded by 
$$\max\{r(S)+1, 2r(S)-1\} \leq \max\{ 2,\ 2(\deg K[S] - \codim K[S])-1 \}.$$
\end{Theorem}

\begin{pf} Let $s= \max\{r(S)+1, 2r(S)-1\}$ and set
$$G =\{ \xbf^\mbf\ybf^\nbf - \xbf^\pbf\ybf^\qbf \in I_\Acal |\  \deg (\xbf^\mbf\ybf^\nbf ) =  
\deg( \xbf^\pbf\ybf^\qbf ) \leq s\}.$$
By \cite[Theorem 1.1]{HS}, $r(S) \le  \deg K[S] - \codim K[S]$.  Hence, it suffices to show that $G$ is a Gr\"obner basis of $ I_\Acal$. In particular this also implies that $G\neq \emptyset$. Assume that this is not the case. Then one can  
find a binomial $b= \xbf^\mbf\ybf^\nbf - \xbf^\pbf\ybf^\qbf\in I_\Acal$ of the smallest  
degree $\deg b>s$ such that $\gin(g) \nmid \xbf^\mbf\ybf^\nbf $ for all $g\in G$.

If $\deg(\xbf^\mbf) \geq r(S)+1$, then we can write $\xbf^\mbf = \xbf^{\mbf'} \xbf^{\mbf"}$,  
where $\deg(\xbf^{\mbf'}) = r(S)+1$. By the definition of $r(S)$ we can find $\mbf^*, \nbf^*$  
such that $\deg(\xbf^{\mbf^*}) =r(S)$ and  $g := \xbf^{\mbf'} - \xbf^{\mbf^*}\ybf^{\nbf^*}  
\in I_\Acal$ (note that  $\xbf^{\mbf'} >  
\xbf^{\mbf^*}\ybf^{\nbf^*}$). Then $g\in G$ and $\gin(g) = \xbf^{\mbf'} \mid  
\xbf^\mbf\ybf^\nbf $, a contradiction. Thus $\deg(\xbf^\mbf) \leq r(S)$. 

If $\deg(\xbf^\pbf)\geq r(S)+1$, then as above, we can find  
$\pbf', \pbf"$ such that $\xbf^\pbf\ybf^\qbf - \xbf^{\pbf'}\ybf^{\pbf" + \qbf} \in I_\Acal$ and  
$\deg(\xbf^{\pbf'}) =r(S) < \deg( \xbf^{\pbf})$. Then
$$\xbf^\mbf\ybf^\nbf - \xbf^{\pbf'}\ybf^{\pbf" + \qbf} = (\xbf^\mbf\ybf^\nbf -  
\xbf^\pbf\ybf^\qbf) + ( \xbf^\pbf\ybf^\qbf - \xbf^{\pbf'}\ybf^{\pbf" + \qbf}) \in I_\Acal,$$ 
and $\xbf^\mbf\ybf^\nbf > \xbf^\pbf\ybf^\qbf > \xbf^{\pbf'}\ybf^{\pbf" + \qbf}$. Hence,  
replacing $\xbf^\pbf\ybf^\qbf $ by $\xbf^{\pbf'}\ybf^{\pbf" + \qbf}$, we may assume from the  
beginning that $\deg(\xbf^\pbf)\leq r(S)$.

Now, since $\xbf^\mbf\ybf^\nbf - \xbf^\pbf\ybf^\qbf\in I_\Acal$, we have
$$ \sum_{i=1}^c m_i\abf_i + \sum_{j=1}^d n_j\ebf_j = \sum_{i=1}^c p_i\abf_i +  
\sum_{j=1}^d q_j\ebf_j. $$
From the minimality of $\deg(\xbf^\mbf\ybf^\nbf)$ we may assume that  
$\xbf^\mbf\ybf^\nbf$ and $\xbf^\pbf\ybf^\qbf$ have no common variable. That means if we  
set  $C= \{ i|\ m_i \neq 0\}$ and $D= \{j|\ n_j \neq 0\}$, then the above equality can be  
rewritten as
$$\sum_{i\in C} m_i\abf_i + \sum_{j\in D} n_j\ebf_j = \sum_{i\not\in C} p_i\abf_i +  
\sum_{j\not\in D} q_j\ebf_j.$$
Hence
$$\sum_{j\in D} \sum_{i\in C} m_ia_{ij} + \sum_{j\in D} n_j\alpha  = \sum_{j\in  
D}\sum_{i\not\in C} p_i a_{ij} = \sum_{i\not\in C} p_i \sum_{j\in D}a_{ij} \leq  
\sum_{i\not\in C} p_i\alpha .$$
This implies 
\begin{equation}\label{EA1}
\sum_{j=1}^d n_j = \sum_{j\in D} n_j \leq \sum_{i\not\in C} p_i = \deg(\xbf^\pbf).
\end{equation}
The equality holds  if  and only if $m_ia_{ij}=0$ for all $(i,j)\in C\times D$ and   
$p_ia_{ij}=0$ for all $(i,j)$ such that $ i\not\in C$ and $ j\not\in D$. This yields $\sum_{i\in  
C} m_i\abf_i = \sum_{j\not\in D} q_j\ebf_j$, which means $\xbf^\mbf - \ybf^\qbf \in  
I_\Acal$. Since $\xbf^\mbf > \ybf^\qbf$ and $\deg(\xbf^\mbf )\leq r(S)$, $g:= \xbf^\mbf -  
\ybf^\qbf \in G$. But this is impossible because $\gin(g) \mid \xbf^\mbf \ybf^\nbf $.
Hence, by (\ref{EA1}), we must have $\sum_{j=1}^d n_j < \deg(\xbf^\pbf)\leq  r(S)$, and so 
$$\deg(b) = \deg( \xbf^\mbf) + \sum_{j=1}^d n_j  \leq 2r(S)-1 \leq s,$$
a contradiction. The theorem is proved. \hfill $\square$
\end{pf}

 It should be noted that if $S$ is not necessarily a simplicial semigroup, then Sturmfels  
\cite{St2} showed that w.r.t. any term order, the maximum degree in a minimal  Gr\"obner basis of  $I_\Acal$ is bounded by $c\cdot \deg K[S]$.

The following example shows that estimations in Theorem \ref{A1} do not hold for an arbitrary term order.  

 \begin{Example} \label{A1A3}{\rm Let $\Acal =\{(4,0),(3,1),(1,3),(0,4)\}$. Then
 $$I_\Acal = ( x_1x_2-y_1y_2,x_1^3-x_2y_1^2,x_2^3 -x_1y_2^2, x_2^2y_1 - x_1^2  
y_2).$$
 This  is also a minimal Gr\"obner basis of $I_\Acal$ w.r.t. the reverse  
lexicographic order. W.r.t. the lexicographic order we get the following minimal Gr\"obner  
basis:
 $$\{ x_1x_2-y_1y_2,x_1^3-x_2y_1^2, x_1y_2^2- x_2^3 ,  x_1^2 y_2 - x_2^2y_1, x_2^4  
-y_1y_2^3\}.$$
 In this example $r(S) = \deg K[S]- \codim K[S] =2$ and both bounds in Theorem \ref{A1}  are equal to $3$.}
 \end{Example}
 
 The above example also provides a case when upper bounds in Theorems \ref{A1} are tight. However if $\deg K[S] - \codim K[S] \geq 3$  we believe that the second bound  is never attained (see the conjecture  mentioned in the introduction). Similarly, we don't think that the first bound  is attained if $r(S)$ is big. However,  
 the following example shows that  in general it is at most twice of the best possible bound.

\begin{Example} \label{A1b} {\rm Given $d\geq 2$ and $\alpha \geq \max\{4, d+1\}$. Let
\begin{equation} \label{EA1b}
\Acal = M_{\alpha ,d}\setminus \{ (\beta , \alpha -\beta ,0,...,0)|\ 2\leq \beta \leq \alpha -2\}.
\end{equation}
We may assume $\abf_1 = (\alpha -1,1,0,...,0)$ and $\abf_2 = (1,\alpha -1,0,...,0)$. If  $S \ni (\alpha -2)\abf_1 =  
\sum m_i\abf_i + \sum n_j\ebf_j$ with $\sum n_j >0$, comparing $d-1$ last coordinates, one  
should have $\alpha -2 = m_1 + m_2 (\alpha -1) + n_2\alpha $. This implies $n_2=m_2 =0$  
and 
$m_1=\alpha -2$, which is impossible, since 
$m_1= (\alpha -2) - \sum n_j <\alpha -2$. Hence $(\alpha -2)\abf_1\not\in \{\ebf_1,...,\ebf_d\}  
+ (\alpha -3)\Acal$ and $r(S) \geq \alpha -2$.

Let $\bbf = (b_1,...,b_d)\in \Nset^d$ such that $\alpha \mid b_1 + \cdots + b_d$ and  
$b_3+\cdots +b_d>0$. By induction on $\deg(\bbf) : = (b_1+\cdots +b_d)/\alpha $, we show  
that $\bbf \in S$. The case $\deg(\bbf) =1$ follows from (\ref{EA1b}). Let $\deg(\bbf) \geq 2$.  
If $b_1\geq \alpha $, then $\bbf = \ebf_1 + \bbf'$ with $b'_3 + \cdots + b'_d >0$.  By the  
induction hypothesis, $\bbf'\in S$ and hence $\bbf\in S$. The same holds if $b_2 \geq  
\alpha$. Hence we may assume that  $b_1,b_2 <\alpha $. In this case $b_2+b_3+\cdots + b_d  
\geq \alpha +1$, and we can find $b'_2= b_2, b'_3 \leq b_3,...,b'_d \leq b_d$ such that  
$b'_2+\cdots + b'_d = \alpha $. Let $b'_1=0$. Then both elements $\bbf'$ and $\bbf-\bbf'$  
satisfy the induction hypothesis, which implies $\bbf = \bbf'+ (\bbf-\bbf') \in S$.

Further, let $\bbf= (b_1,b_2,0,...,0)$ with $b_1+b_2 = \alpha (\alpha -2)$.  We  show that  also
$\bbf\in S$. Indeed,  we can write $b_2 = p\alpha +q$, where $p\leq \alpha -2, \ q\leq  \alpha  
-1$. Note that $p= \alpha -2$ implies $q=0$ and $\bbf = (\alpha -2)\ebf_2 \in S$. 

Let $p\le  \alpha -3$. If $p+q \geq  \alpha -1$, then
$$\bbf = 0\abf_1 + (\alpha -q) \abf_2 + (\alpha -3-p)\ebf_1 + (p+q-\alpha +1)\ebf_2\in S.$$
Otherwise ($p+q \leq  \alpha -2$),
$$\bbf = q\abf_1 + 0 \abf_2 + (\alpha - 2-p-q)\ebf_1 + p\ebf_2\in S.$$
 Summarizing the above arguments we get that $\bbf \in S$ if $\deg (\bbf) =\alpha -2$. 

Now let $\abf \in (\alpha -1)\Acal$. Since $\alpha \geq d+1$,   $a_1+ \cdots + a_d = \alpha  
(\alpha -1) \geq d\alpha $ and there is an index $i$ such that $a_i \geq \alpha $. Note that  
$\deg (\abf - \ebf_i) = \alpha -2$. By the above result $\abf - \ebf_i \in S$. Hence $\abf = \ebf_i  
+ (\abf-\ebf_i) \in \{\ebf_1,...,\ebf_d\} + S$, which implies $r(S) \leq \alpha -2$.

Summing up we get $r(S) = \alpha -2$. 

On the other hand,   $x_1^{\alpha -1} -  
x_2y_1^{\alpha -2} \in I_\Acal$,  $x_1^{\alpha -1} > x_2y_1^{\alpha -2}$ and there is no  
other binomial of $I_\Acal$ whose first term divides $x_1^{\alpha -1}$. Therefore the  
binomial $x_1^{\alpha -1} - x_2y_1^{\alpha -2}$ must be contained in the reduced  
Gr\"obner basis of $I_\Acal$. The degree of this binomial is $\alpha -1 $, while  the  first bound of Theorem \ref{A1} is $2\alpha -5$. 
Note that the Eisenbud-Goto bound in this example is $\alpha^{d-1} + \alpha +d -{\alpha +d-1\choose d-1} -2$. 
}\end{Example}

 It was also shown  that the Castelnuovo-Mumford regularity of  
$\reg(I_\Acal)$ is bounded by $c(\alpha -1) +1$ (see \cite[Theorem 3.2(i)]{HS}).  In the following theorem we obtain a similar result for Gr\"obner bases.
 
 \begin{Theorem} \label{A4} The maximum degree in a minimal  Gr\"obner basis of $ I_\Acal$ is bounded by  $\max\{c, \alpha ,\ c(\alpha -1) -1\} \leq c(\alpha -1)$.
 \end{Theorem}
 
 \begin{pf} The proof is similar to that of Theorem \ref{A1}. Let $s= \max\{c, \alpha ,\ c(\alpha  
-1) -1\}$ and set
$$G =\{ \xbf^\mbf\ybf^\nbf - \xbf^\pbf\ybf^\qbf \in I_\Acal |\  \deg (\xbf^\mbf\ybf^\nbf ) =  
\deg( \xbf^\pbf\ybf^\qbf ) \leq s\}.$$
 Assume that $G$ is not a Gr\"obner basis. Then one can find a binomial $b=  
\xbf^\mbf\ybf^\nbf - \xbf^\pbf\ybf^\qbf\in I_\Acal$ of the smallest degree $\deg b>s$ such  
that $\gin(g) \nmid \xbf^\mbf\ybf^\nbf $ for all $g\in G$. Since $\alpha \abf_i = a_{i1}\ebf_1  
+ \cdots + a_{id}\ebf_d$, $x_i^\alpha - \ybf^{\abf_i} \in G$ for all $i=1,...,c$. Note that  
$x_i^\alpha > \ybf^{\abf_i}$. Since $\gin(x_i^\alpha - \ybf^{\abf_i}) \nmid  
\xbf^\mbf\ybf^\nbf $, we must have  $m_i \leq \alpha -1$ for all $i\leq c$.  

If $p_i \geq \alpha $,  then
$$\xbf^\mbf\ybf^\nbf  - \frac{\xbf^\pbf}{x_i^\alpha}\ybf^{\qbf + \abf_i} =  
(\xbf^\mbf\ybf^\nbf - \xbf^\pbf\ybf^\qbf) + (x_i^\alpha - \ybf^{\abf_i})  
\frac{\xbf^\pbf}{x_i^\alpha}\ybf^{\qbf} \in I_{\Acal}.$$
Note that $\xbf^\mbf\ybf^\nbf > \xbf^\pbf\ybf^\qbf > \frac{\xbf^\pbf}{x_i^\alpha}\ybf^{\qbf  
+ \abf_i}$. Replacing $b$ by $\xbf^\mbf\ybf^\nbf  - \frac{\xbf^\pbf}{x_i^\alpha}\ybf^{\qbf +  
\abf_i}$ and repeating this procedure, we may also assume that $p_i \leq \alpha -1$ for all  
$i\leq c$.  

As in the proof of Theorem \ref{A1}, let $C= \{ i|\ m_i \neq 0\}$ and $D= \{j|\ n_j \neq 0\}$.  
Then  we can also conclude that 
\begin{equation} \label{EA4}
\sum_{j\in D} n_j \leq \sum_{i\not\in C} p_i \leq (c - \sharp C) (\alpha -1),
\end{equation}
and that $\sum_{j\in D} n_j  = (c - \sharp C) (\alpha -1)$ implies  $\xbf^\mbf - \ybf^\qbf \in  
I_\Acal$.  Hence
$$\deg (\xbf^\mbf\ybf^\nbf ) = \sum_{i\in C}m_i + \sum_{j\in D} n_j \leq  \sharp C (\alpha  
-1) + (c - \sharp C) (\alpha -1) = c(\alpha -1).$$
Since $\deg (\xbf^\mbf\ybf^\nbf ) = \deg(b) \geq c(\alpha -1)$, we must have $\deg  
(\xbf^\mbf\ybf^\nbf ) = c(\alpha -1)$. Therefore $\sum_{j\in D} n_j  = (c - \sharp C) (\alpha  
-1)$ and $m_i  = \alpha -1$ for all $i\in C$.  By (\ref{EA4}) we have  $\xbf^\mbf - \ybf^\qbf  
\in I_\Acal$. If $C \neq \{1,...,c\}$, then $\deg(\xbf^\mbf ) \leq s$ and $\xbf^\mbf - \ybf^\qbf  
\in G$, which is impossible because $\xbf^\mbf \mid \gin(b)$. Thus $C= \{1,...,c\}$. This  
yields $D= \emptyset $ and
$$b = (x_1\cdots x_c)^{\alpha -1} - \ybf^\qbf.$$
Let $\abf = \abf_1 + \cdots + \abf_c$ and $\abf := (a_1,...,a_d)$. The above equality assures that $\alpha \mid (\alpha  
-1)a_i $ for all $i=1,...,c$. This implies $a_i = q'_i\alpha $ for some $q'_i\in \Nset$. But then  
$g:= x_1\cdots x_c - y_1^{q'_1}\cdots y_d^{q'_d} \in I_\Acal$. Since $\deg(x_1\cdots x_c)  
= c \leq s$, $g\in G$ and we get a contradiction that $\gin(g) =  x_1\cdots x_c\mid \gin(b) = (  
x_1\cdots x_c)^{\alpha -1}$. The proof of the theorem is completed. \hfill $\square$
\end{pf}

\begin{Remark} {\rm  Theorem \ref{A4} shows that our conjecture holds if the codimension is not too big. This includes the case $c=2$ and  $\deg K[S]  > \alpha $, because we always have   $\alpha \mid \deg K[S]$ (see \cite[Lemma 3.4]{HS}).  Note that the case $c=2$  (even if $\deg K[S]  = \alpha $) was completely solved by Peeva and Sturmfels   (see \cite[Theorem  
7.3 and Proposition 8.3]{PS}). Another proof was recently given in \cite{BGM} (see Theorems  
2.1, 2.8 and 3.5 there).
}\end{Remark}

 Note that  an ideal is usually given by its generating set and this  set  serves as the input data for computing a Gr\"obner basis. In the case of a toric ideal $I_\Acal$ the input data is  $\Acal$, and before computing a Gr\"obner basis of $I_\Acal$ we have to compute  a generating set of this ideal. However, in many algorithms we get a Gr\"obner basis of $I_\Acal$ as a by-product of  computing a generating set. The last result of this section shows that, by using a suitable term order, the computation of  
simplicial toric ideals runs rather quickly.
  In order to compute $I_\Acal$, a standard procedure is the following (see, e.g., \cite{St1},  
Algorithm 4.5):
  
  \begin{itemize}
  \item[1.] Form the ideal $J_\Acal = (x_1-\tbf^{\abf_1},...,x_c-\tbf^{\abf_c}, y_1-t_1^\alpha  
,...,y_d-t_d^\alpha ) \subset K[\tbf,\xbf,\ybf]$.
  \item[2.] Compute a Gr\"obner basis $G'$ of $J_\Acal$ by Buchberger's algorithm, using an  
elimination order $\preceq $ with respect to the variables $t_1,...,t_d$. Here we assume $t_1\succeq  \cdots \succeq t_d \succeq  x_1\succeq  \cdots \succeq  y_d$.
  \item[3.] From $G'$ get a Gr\"obner basis $G= G' \cap K[\xbf,\ybf]$ of $I_\Acal = J_\Acal  
\cap K[\xbf,\ybf]$.
  \end{itemize}
  
 Though this algorithm is  not the best one, the following result says that it requires not too many steps. Moreover, by Theorems \ref{A1} and  \ref{A4}, in order to compute $G$ it is sufficient  to compute those elements of  $G'$ which have degrees up to $\min\{2r(S) , c(\alpha -1)\}$, that means we can do truncation in the above algorithm.
  
  \begin{Proposition}\label{A6} Assume that the restriction of the elimination order $\preceq  
$ on $K[\xbf,\ybf]$  is the reverse lexicographic order. Then  the maximum degree in a minimal  Gr\"obner basis of
$J_\Acal$ is bounded by $d(\alpha -1) + \min\{2r(S) , c(\alpha -1)\}$.
  \end{Proposition}
  
  \begin{pf} The proof is similar to that of Theorem \ref{A1}.  We give here a sketch. Let $s= d(\alpha -1) + \min\{2r(S)  
, c(\alpha -1)\}$ and set
$$G =\{ \tbf^\pbf\xbf^\mbf\ybf^\nbf - \tbf^{\pbf'}\xbf^{\mbf'}\ybf^{\nbf'} \in J_\Acal |\  \deg  
(\tbf^\pbf\xbf^\mbf\ybf^\nbf - \tbf^{\pbf'}\xbf^{\mbf'}\ybf^{\nbf'} ) \leq s\}.$$
 Assume that $G$ is not a Gr\"obner basis. Then one can find a binomial $b=  
\tbf^\pbf\xbf^\mbf\ybf^\nbf - \tbf^{\pbf'}\xbf^{\mbf'}\ybf^{\nbf'}\in J_\Acal$ of the smallest  
degree $\deg b>s$ such that 
$\gin(g) \nmid \tbf^\pbf \xbf^\mbf\ybf^\nbf $ for all $g\in G$. Since $t_i^\alpha - y_i\in G$  
and $t_i^\alpha \succ  y_i$,  as in the proof of Theorem \ref{A4}, we can assume that  $p_i, p'_i\leq \alpha -1$ for all $i\leq d$. 

Using arguments in the proof of Theorem \ref{A1} we may assume that $\deg(\xbf^\mbf),$ $ \deg(\xbf^{\mbf'}) \leq r(S)$. Note that $J_\Acal$ is the kernel of the epimorphism $K[\tbf,\xbf,\ybf] \rightarrow K[\tbf]$ mapping $t_j, x_i, y_j$ to $t_j,  \tbf^{\abf_i}, t_j^\alpha$, respectively. Therefore $b\in J_\Acal$ if and only if
\begin{equation} \label{EA6}  p_j +  \sum_{i\in C} m_i a_{ij} +  n_j\alpha =  p'_j +  \sum_{i\not\in C} m'_i a_{ij} +  n'_j\alpha
\end{equation}
for all $j \le d$, where $C$ and $D$ are the same as in the proof of Theorem \ref{A1}. Then, instead of (\ref{EA1}) we get
$$ \sum_{j\in D} p_j + \sum_{j=1}^d n_j \leq   \sum_{j\in D} p'_j  + \deg(\xbf^{\mbf'}) \le (\sharp D)(\alpha -1) +  r(S),$$
which implies
$$\sum_{j=1}^d p_j + \sum_{j=1}^d n_j \leq   d(\alpha -1) +  r(S).$$
Hence $\deg(\tbf^\pbf\xbf^\mbf\ybf^\nbf ) \le d(\alpha -1) +  2r(S)$. Similarly, $\deg (\tbf^{\pbf'}\xbf^{\mbf'}\ybf^{\nbf'}) \leq d(\alpha -1) +  2r(S)$, and so $\deg(b) \le d(\alpha -1) +  2r(S)$.

Now, applying arguments in the proof of Theorem \ref{A4} to (\ref{EA6})  we can also conclude that
$\deg (b) \leq d(\alpha -1) +  c(\alpha  
-1)$. 

Summing up, we get $\deg(b) \leq s$,  a contradiction. \hfill $\square$
 \end{pf}

\section{Eisenbud-Goto bound}\label{B}

In this section we will provide some partial positive answers to our conjecture.

Recall that a quotient ring $R/I$ modulo a homogeneous  ideal $I$ is said to be a {\it  
generalized Cohen-Macaulay} ring if all local cohomology modules $H^i_\mfr(R/I),\ i<\dim  
R/I,$ with the support in the maximal homogeneous ideal $\mfr$ of $R/I$ are of finite length  
(see the Appendix in \cite{SV1}). The Castelnuovo-Mumford regularity of a finitely generated  
graded $R$-module $M$ is the number 
$$\reg(M) = \max\{ n|\ [H^i_\mfr(M)]_{n-i} \neq 0 \ \text{for} \ i\geq 0\}.$$
Note that $\reg(I) = \reg(R/I) + 1$. The following result is a simple observation, but has some  
interesting consequences.

\begin{Lemma} \label{C1} Assume that $K[S]$ is a generalized Cohen-Macaulay ring. Then  
the maximum degree in a minimal  Gr\"obner basis of $I_\Acal$ is bounded by $\reg  
I_\Acal$.
\end{Lemma}

\begin{pf} Note that $y_1,...,y_d$ is a system of parameters of $K[S]$. Since $K[S] \cong   
K[\xbf, \ybf]/I_\Acal$ is a generalized Cohen-Macaulay ring, the ideal $I_\Acal$ and all ideals  
$(I_\Acal, y_d,...,y_i),\ i=d,d-1,...,1$, are unmixed  up to  $\mfr$-primary components (see  
\cite[Proposition 3 in the Appendix]{SV1}). In particular, $y_{i-1}$ is a non-zero divisor on the  
ring $K[\xbf, \ybf]/(I_\Acal, y_d,...,y_i)^{\rm sat}$, where $J^{\rm sat}= \cup_{n\geq 1} J:\mfr^n$  
denotes the saturation of $J$. This means $y_d,...,y_1$ is a generic sequence of $K[S]$ in the  
sense of \cite[Definition 1.5]{BS}. By \cite[Corollary 2.5]{BS},  the maximum degree in  a  
minimal Gr\"obner basis of $I_\Acal$ is bounded by $\reg(I_\Acal)$. \hfill $\square$
\end{pf}

\begin{Remark}\label{C1b} {\rm  Note that $y_d,...,y_1$ is  always a system of  
parameters of $\gin (I_\Acal)$. This follows from the fact that $x_i^\alpha \in \gin (I_\Acal)$  
for all $i\leq c$ (since $x_i^\alpha - \ybf^{\abf_i} \in I_\Acal$). However, if   $K[S]$ is not a  
generalized Cohen-Macaulay ring, it maybe no more  a generic sequence of  
$K[S]$. For example, let $d=\alpha =3$ and
$$\begin{array}{ll}  \Acal = &\{\ebf_1,\ebf_2,\ebf_3,\ \abf_1= (2,0,1),\ \abf_2 = (1,2,0),  \
\abf_3 = (1,1,1),\\ &\quad \abf_4=(1,0,2),\ \abf_5=(0,2,1),\ \abf_6=(0,1,2)\}.\end{array}$$
Then 
$$\begin{array}{ll} \gin(I_\Acal) = & (x_1x_2, x_2x_3, x_2x_5, x_1^2, x_1x_3, x_3^2,  
x_2x_4, x_2x_6, x_3x_5, x_5^2, x_1x_4, x_3x_4,\\ & \quad x_4x_5, x_4^2, x_3x_6, x_5x_6,  
x_4x_6, x_6^2, x_2^3, x_1x_6y_2) .
\end{array}$$
 Clearly $(\gin(I_\Acal),y_3)^{sat} = (\gin(I_\Acal),y_3)$ and $y_2$ is a zero divisor of  
$K[\xbf]/(\gin(I_\Acal), y_3)$. Hence, by \cite[Theorem 2.4(a)]{BS}, $y_3,y_2,y_1$ is not a  
generic sequence of  $K[\xbf]/I_\Acal$.}
\end{Remark}

\begin{Corollary} \label{C3} The maximum degree in a minimal  Gr\"obner basis of $ I_\Acal$ is bounded by $\deg K[S] - \codim K[S]+1$ in the following cases:
\begin{itemize}
\item[(i)] $d=2$,
\item[(ii)]  $K[S]$ is a so-called Buchsbaum ring,
\item[(iii)]  $K[S]$ is  a simplicial semigroup ring with isolated  
singularity, or equivalently,  $\Acal$ contains all points of  
$M_{\alpha ,d}$ of type $(0,..,\alpha -1,...,1,...,0)$, where $\alpha -1,1$ stay in the $i$-th and  
$j$-th positions, respectively, and the other coordinates are zero.
\end{itemize}
\end{Corollary}

\begin{pf} In all these cases,  $K[S]$ is a generalized Cohen-Macaulay ring and it is known that $\reg(I_\Acal) \leq \deg K[S] - \codim K[S]+1$ (the case (i)  is proved  in \cite{GLP},  (ii) in  \cite[Theorem 1]{SV2} and (iii) is  \cite[ Corollary 2.2]{HH}). Hence the statement follows from Lemma \ref{C1}. \hfill $\square$
\end{pf}

Recall that a semigroup $S$ is said to be {\it normal} if $S = \Zset(S) \cap \Nset^d$. Under this  
condition, it is well-known  that $\reg(I_\Acal) \leq d$ (this holds even without the  
assumption $S$ being simplicial, see \cite[Proposition 13.14]{St1}).  Hence, by Lemma \ref{C1},  the maximum degree in a minimal  Gr\"obner basis of $I_\Acal$ is bounded by $d$. This gives  a partial answer to  the following question posed by Sturmfels in \cite[p. 136]{St1}: If the semigroup $S$ is normal, does the toric ideal  $I_\Acal$ posses a Gr\"obner basis of degree at most $d$?

Under the assumption of the following result it was shown in \cite[Proposition 3.7]{HS} that  
$\reg(I_\Acal) \leq \deg K[S] - \codim K[S]+1$. Unfortunately we cannot use it to derive the  
corresponding result for a Gr\"obner basis, because $K[S]$ is in general not a generalized  
Cohen-Macaulay ring.

\begin{Proposition} \label{B2} Assume that $\deg K[S] = \alpha^{d-1}$ and $\alpha \leq  
d-1$. Then the maximum degree in a minimal  Gr\"obner basis of $I_\Acal$ is bounded by $\deg K[S] - \codim K[S]+1$.
\end{Proposition}

In order to prove this proposition we need to recall a result from \cite{HS}. Let $\Pcal$ denote  
the convex polytope spanned by $\Acal \subset \Rset^d$. Note that $\Pcal$ is a
$(d-1)$-dimensional polytope whose faces are spanned by 
$$\Acal_I = \{\abf\in \Acal|\ a_{i} = 0\ \ \text{for\ all\ } i\in I\},$$
where $I \subseteq \{1,...,d\}$. Let
 $\Pcal_I$ denote the corresponding  face of $\Pcal$. (For short, we will also write $\Acal_i,\  
\Pcal_i$ instead of $\Acal_{\{i\}},\ \Pcal_{\{i\}}$.)
  We  say that a face $\Pcal_I$ is {\it full} if $\Acal_I$ contains all points of 
$M_{\alpha, d}$ lying on this face, i.e. if $\Acal_I = \Pcal_I \cap M_{\alpha,d}$. 

\begin{Lemma} {\rm (\cite[ Lemma 1.2]{HS})}. \label{A2}  If $\Pcal$ has a full face of dimension $i$, then
$r(S) \leq \alpha^{d-1-i} + i-1.$
 \end{Lemma}
 
\noindent PROOF OF PROPOSITION \ref{B2}.  If $\Acal = M_{\alpha ,d}$, then by Corollary \ref{C3}(iii) we are  
done.  Hence we may assume that
$$c\leq \sharp M_{\alpha,d}-1-d = {\alpha +d-1\choose d-1}-d-1.$$
If $\alpha \geq 3,\ d\geq 6$ or $\alpha =4,\ d=5$, then by \cite[Claim 1, p. 141]{HS},  $c\leq  
\alpha^{d-2}$. Hence
$$\deg K[S] - c+1 > \alpha^{d-1} - \alpha^{d-2} = \alpha^{d-2}(\alpha -1) \geq c(\alpha -1),$$
and by Theorem \ref{A4} we are done. Thus the left cases are: $\alpha =2, \ d\geq 3$; $\alpha  
=3,\ d=4$ and $\alpha =3,\ d=5$. We consider these cases separately.
\vskip0.3cm

{\bf Case 1}: $\alpha =2, \ d\geq 3$. Then $c\leq d(d+1)/2 - (d+1) =(d-2)(d+1)/2$. It is easy to  
verify that $(d-2)(d+1)/2-1\leq 2^{d-2}$. Hence, if $c\leq (d-2)(d+1)/2-1$ we have $c\leq  
2^{d-1} - c = \deg K[S] -c$. By Theorem \ref{A4} we are done. The left case is  
$c=(d-2)(d+1)/2$, i.e. $\Acal$ is obtained from $M_{2,d}$ by deleting exactly one point. We  
may assume $\Acal = M_{2,d}\setminus \{\bbf := (1,1,0,...,0)\}$. Note that $2\abf_i \in  
\{\ebf_1,...,\ebf_d\} + \{\ebf_1,...,\ebf_d\} $ for all $i\leq c$. Moreover, if $\abf_i,\abf_j$ are  
two different points and  $\abf_i,\abf_j,\bbf$ do not lie in the same 2-dimensional face of  
$\Pcal$, then $\abf_i + \abf_j \in \Acal + \{\ebf_1,...,\ebf_d\} $ (see Fig. 1). From this it  
follows that $r(S)=2$. By Theorem \ref{A1}, $I_\Acal$ has a 
Gr\"obner basis of degree at most $3 \leq 2^{d-1} - (d-2)(d+1)/2 + 1 = \deg K[S] -c +1$. 

\centerline{\setlength{\unitlength}{0.4cm}
 \begin{picture}(8,7)
\put(0,2){\line(3,4){3}}
\put(-0.8,2){$\ebf_2$}
\put(1.5,4){\circle{0.2}} \put(1,4){$\bbf$}
\put(0,2){\line(3,-2){3}}
\put(3.1,6){$\ebf_1$}
\multiput(0,2)(0.3,-0.05){20}{.}
\put(3,0){\line(0,1){6}}
\put(2.8,-0.5){$\ebf_3$}
\put(3,0){\line(3,1){3}}
\put(3,6){\line(3,-5){3}}
\put(6.1,1){$\ebf_4$}
\put(3,3){\circle*{0.2}} \put(3.2,3){$\abf_k$}
\put(4.5,0.5){\circle*{0.2}} \put(4.5,-0.1){$\abf_i$}
\put(4.5,3.5){\circle*{0.2}} \put(4.7,3.5){$\abf_j$}
\put(1,-2){ Fig. 1}
\end{picture} 
\begin{picture}(11,7)
\put(3,0){\line(1,0){6}}
\put(3,0){\line(1,2){3}}
\put(6,6){\line(1,-2){3}}
\put(6,6){$\ebf_1$} \put(9,0){$\ebf_3$} \put(2.2,0){$\ebf_2$}
\put(5,0){\circle{0.2}} 
\put(7,0){\circle*{0.2}} \put(6.8,-0.6){$\abf_3$}
\put(4,2){\circle{0.2}} 
\put(5,4){\circle*{0.2}} \put(4,4){$\abf_1$}
\put(8,2){\circle*{0.2}} \put(8,2){$\abf_2$}
\put(7,4){\circle{0.2}} 
\put(6,2){\circle*{0.2}} \put(6.2,2){$\bbf_4$}
\put(4,-2){ Fig. 2}
\end{picture}}
\vskip1cm  

{\bf Case 2}: $\alpha = 3, d=4$. Then $\deg K[S] = 27$ and $c\leq 15$. By Theorem \ref{A4},  
the statement of the proposition holds true for $c\leq 9$. Let $c\geq 10$, i.e. $\Acal$ is  
obtained from $M_{3,4}$ by deleting at most 6 points. We distinguish two subcases.
\vskip0.3cm

{\it Subcase 2a}: Each edge of $\Pcal$ contains exactly one deleting point. In this case  
$c=10$. By Theorem \ref{A1} it suffices to show that $r(S) \leq 8$.

Consider, for example, the facet $\Pcal_4 = \{\abf \in \Pcal|\ a_4 =0\}$. Then $\Acal_4 = \Acal  
\cap \Pcal_4$ has exactly 7 points, say
$\Acal_4 = \{\ebf_1, \ebf_2,\ebf_3, \abf_1,\abf_2,\abf_3,\bbf_4: = (1,1,1,0)\}$ as shown in  
Fig. 2, where $\abf_1$ can be taken as $(2,1,0,0)$, while there are two choices for each of  
$\abf_2$ and $\abf_3$. One can check by computer that in this case the reduction number  
$r(\langle \Acal_4\rangle) \leq 3$. In particular,  $\sum_{i=1}^3m_i\abf_i  +  
n_4\bbf_4\not\in  S + \{\ebf_1,\ebf_2, \ebf_3\}$ implies that $m_i,n_4\leq 2$ and
\begin{equation}\label{EB2}
\sum_{i=1}^3m_i + n_4 \leq 3.
\end{equation}
 Moreover, since
$2\bbf_4 + \abf_1 = (4,3,2,0)  = \ebf_2 + 2(2,0,1,0) = \ebf_1 + \ebf_2 + (1,0,2,0)$ and one of  
two points $(2,0,1,0) $ and $(1,0,2,0)$ on the edge $\overline{\ebf_1\ebf_3}$ must belong to  
$\Acal_4$, we get that $2\bbf_4 + \abf_1 \in S + \{\ebf_1,\ebf_2, \ebf_3\}$. The same is true  
for $2\bbf_4 +\abf_2$ and $2\bbf_4 +\abf_3$. This means, in addition to (\ref{EB2}) we also  
have $m_1=m_2=m_3 =0$ if $n_4=2$.

Finally, we can write $\Acal = \{\ebf_1, \ebf_2,\ebf_3, \abf_1,...,\abf_6,\bbf_1,...,\bbf_4 \},$  
where $\bbf_i$ is the inner point of the facet $\Pcal_i$. Assume that 
$$\sum_{i=1}^6m_i\abf_i  + \sum_{j=1}^4 n_j\bbf_j \not\in  S +  
\{\ebf_1,\ebf_2,\ebf_3,\ebf_4 \}.$$
Then inequalities of Type (\ref{EB2}) should hold for all  facets of $\Pcal$.  If  $\sharp \{j|\ 1\le j\le 4,\ n_j \leq 1\} \ge 3$,  adding all of them we get
$$2\sum_{i=1}^6m_i +  \sum_{j=1}^4 n_j  + \sum_{j=1}^4 n_j  \leq 12 + 1+1+1+2 = 17.$$
Hence $\sum_{i=1}^6m_i +  \sum_{j=1}^4 n_j \leq 8$. 

Otherwise, we may assume either (a) $n_1,n_2 \leq 1$ and $n_3=n_4 = 2$ or (b) $n_2=n_3=n_4 = 2$ holds. In the case (a), we may further assume that 
$m_1 = \cdots = m_5 =0$ and $m_6\le 2$. Then
$$\sum_{i=1}^6m_i +  \sum_{j=1}^4 n_j \leq  m_6 + 1+1+2+2 \leq 8.$$
In the case (b), we have  $m_1 = \cdots = m_6 =0$. Hence
$$\sum_{i=1}^6m_i +   \sum_{j=1}^4 n_j \leq 2+2+2+2 =8.$$
Thus, in all cases we get that $\sum_{i=1}^6m_i +  \sum_{j=1}^4 n_j \leq 8$, which implies $r(S) \leq 8$.
\vskip0.3cm

{\it Subcase 2b}: At least one edge of $\Pcal$ is full. By Lemma \ref{A2}, $r(S)\leq 9$.  
Hence, by Theorem \ref{A1}, the statement holds true if $c\leq 11$.  Moreover, if $\Pcal$ has  
a full facet, then again by Lemma \ref{A2}, $r(S)\leq 5$, and by Theorem \ref{A1} we are  
done. Hence, we may assume that $c=12,13,14$, and $\Pcal$ has no full facet. This  
corresponds to the situation when $\Pcal$ has 2,3 or 4 deleting points.

Assume that $\Pcal$ has a facet, say $\Pcal_4$, which contains exactly one deleting point  
$\bbf$, i.e. one can write $\Acal_4 = \{\ebf_1, \ebf_2,\ebf_3,\abf_1,...,\abf_6\}$ and  
$\bbf\not\in \Acal_4$. By Theorem \ref{A1}, it suffices to show that $r(S)\leq 5$. If this is not  
the case, then one can find $m_1,...,m_c\in \Nset$ such that $\sum_{i=1}^c m_i=6$ and 
\begin{equation} \label{EB2b}
\sum_{i=1}^c m_i\abf_i \not\in S+ \{\ebf_1,\ebf_2,\ebf_3,\ebf_4\}.
\end{equation}
We follow the idea in the proof of \cite[Lemma 1.2]{HS}. Considering $6$  partial sums  
$\abf_1,...,m_1\abf_1, \ m_1\abf_1+\abf_2,...,\ \sum_{i=1}^c m_i\abf_i$ we can find either  
two partial sums whose last coordinates are divisible by $3$, or three partial sums whose last  
coordinates are congruent modulo $3$. Taking also the differences of these partial sums, we can  
find in both cases two partial sums $\bbf_1 = \sum p_i\abf_i$ and $\bbf_2 = \sum q_i \abf_i$  
such that $m_i \ge p_i \ge q_i \ge 0,\ (i\leq 3)$, $\deg(\bbf_1) > \deg(\bbf_2) \geq 2$ and the last  
coordinates of $\bbf_1,\bbf_2$ are divisible by $3$.  Note that $b_3: = b_1-b_2 \neq 0$  also is a partial sum of $\sum_{i=1}^c m_i\abf_i$.  Fix $i\in \{2,3\}$. We can write $\bbf_i = \bbf_i' +  n_i \ebf_4$, where $\bbf_i' \in \langle \Acal_4, \bbf\rangle$. By (\ref{EB2b}) we  
must have $\bbf_i \not\in S+ \{\ebf_1,\ebf_2,\ebf_3,\ebf_4\}$, which yields $0 \neq  
\bbf_i'\not\in \{\ebf_1,\ebf_2,\ebf_3\} + \langle \Acal_4\rangle$. Together with the  
fact  $2\bbf \in \langle \Acal_4\rangle$, this implies $\bbf_i' \in \{\bbf, \abf_1,...,\abf_6\} +  
\langle \Acal_4\rangle$. Since also all elements $2\bbf, \bbf+\abf_1,...,\bbf+\abf_6 \in \langle  
\Acal_4\rangle$, the previous relation assures that $\bbf_2'+\bbf_3'\in \langle \Acal_4\rangle  
\subset S$. By (\ref{EB2b}) we must have $n_2=n_3 =0$, and so $\bbf_1= \bbf_2+\bbf_3= \bbf_2'+\bbf_3' \in \langle  
\Acal_4\rangle$. However it is easy (or using computer) to see  that $r(\langle \Acal_4\rangle)  
=2$. Since $\deg(\bbf_1)\geq 3$, 
$\bbf_1\in \langle \Acal_4\rangle + \{\ebf_1,\ebf_2,\ebf_3\} \subseteq S+  
\{\ebf_1,\ebf_2,\ebf_3,\ebf_4\},$
which contradicts (\ref{EB2b}).

Thus, each facet of $\Pcal$ must have at least two deleting points. In particular, $c=12$ and  
$\Pcal$ has exactly 4 deleting points. There are only two situations shown in Fig. 3 and Fig. 4.  
In the situation of Fig. 3 there is eventually one configuration, and by computer we see that  
$r(S)=2$. In the situation of Fig. 4 one can show as in Subcase 2a (or using computer for eight  
different configurations), that $r(S)\leq 8$. But then by Theorem \ref{A1}, $I_\Acal$ has a  
Gr\"obner basis of degree at most $15 < \deg K[S] - c+ 1= 16$. The Subcase 2b is completely  
solved.
 
\centerline{\setlength{\unitlength}{0.4cm}
 \begin{picture}(8,7)
\put(0,2){\line(3,4){3}}
\put(1,3.3){\circle{0.2}} 
\put(2,4.7){\circle{0.2}}
\put(0,2){\line(3,-2){3}}
\multiput(0,2)(0.3,-0.05){20}{.}
\put(3,0){\line(0,1){6}}
\put(3,0){\line(3,1){3}}
\put(3,6){\line(3,-5){3}}
\put(4,0.3){\circle{0.2}} 
\put(5,0.7){\circle{0.2}} 
\put(1,-2){ Fig. 3}
\end{picture} 
\begin{picture}(8,7)
\put(0,2){\line(3,4){3}}
\put(1,3.3){\circle{0.2}} 
\put(3,4){\circle{0.2}}
\put(0,2){\line(3,-2){3}}
\multiput(0,2)(0.3,-0.05){20}{.}
\put(3,0){\line(0,1){6}}
\put(3,0){\line(3,1){3}}
\put(3,6){\line(3,-5){3}}
\put(4,0.3){\circle{0.2}} 
\put(2,1.7){\circle{0.2}} 
\put(1,-2){ Fig. 4}
\end{picture} }
\vskip1cm

{\bf Case 3}: $\alpha =3, d=5$. We have $c\leq 29$ and $\deg K[S] = 81$. If $c\leq 27$, then  
$\deg K[S] - c + 1 \geq 54 \geq 2c$, and by Theorem \ref{A4} we are done. If $c= 28, 29$,  
then $\Acal$ is obtained from $M_{3,5}$ by deleting 1 or 2 points. But then $\Pcal$ has a  
full $2$-dimensional face. By Lemma \ref{A2}, $r(S) \leq 10$. Hence, by Theorem  
\ref{A1}, we are  also done in this subcase. \hfill $\square$
\vskip0.5cm

Finally we show that if on an edge of $\Pcal$ there are enough points belonging to $\Acal$,  
then the Eisenbud-Goto bound also holds for the maximum degree in a minimal Gr\"obner basis of  
$I_\Acal$. Note that in this setting, the Eisenbud-Goto conjecture on $\reg(I_\Acal)$  is still  
not verified (cf. \cite[Corollary 3.8]{HS}).

\begin{Proposition} \label{B3} Assume that $\deg K[S]= \alpha^{d-1}$ and there exists an  
edge of $\Pcal$ such that it is either full or at least $(\frac{3}{4} + \frac{1}{4d})\alpha +2$  
integer points on it belong to $\Acal$. Then the maximum degree in a minimal  Gr\"obner basis of $I_\Acal$ is bounded by $\deg K[S] - \codim K[S] +1$.
\end{Proposition}

\begin{pf} By Corollary \ref{C3}  and Proposition \ref{B2} we may assume  
that $\alpha \geq d\geq 3$ and
$$c\leq {\alpha + d-1 \choose d-1} - d-1.$$
First we consider the case when at least $(\frac{3}{4} + \frac{1}{4d})\alpha +2$ integer  
points on an edge belong to $\Acal$. By \cite[Lemma 1.3]{HS}, $r(S) \leq  
\frac{d-1}{4d}\alpha^{d-1}$. Hence, by Theorem \ref{A1}, it suffices to show that
$$\alpha^{d-1} - {\alpha + d-1 \choose d-1} + d+1 \geq \frac{d-1}{2d}\alpha^{d-1} -1,$$
or equivalently
\begin{equation}\label{EB3}
\frac{d+1}{2d}\alpha^{d-1} + d+2 \geq {\alpha + d-1 \choose d-1}.
\end{equation}
We show this by induction on $d\geq 3$. For $d=3$ this is equivalent to $\alpha^2 - 9\alpha  
+24\geq 0$. So, assume that the inequality holds for $d\geq 3$. In the dimension $d+1$, by  
induction we have
$$\begin{array}{ll} 
 \displaystyle{{\alpha + d \choose d}} &=  \displaystyle{ \frac{\alpha +d}{d}{\alpha + d-1  
\choose d-1}
 \leq \frac{\alpha +d}{d}(\frac{d+1}{2d}\alpha^{d-1} + d+2) }  \\
&=  \displaystyle{ \frac{(\alpha +d)(d+1)}{2d^2}\alpha^{d-1} + \frac{d+2}{d}\alpha + d+2.}
\end{array}$$
Hence
$$\begin{array}{l} 
\displaystyle{\frac{d+2}{2(d+1)}\alpha^d + d+3 - {\alpha + d \choose d}}\\
\quad \quad  \geq \displaystyle{ \alpha^{d-1}\left[\frac{d+2}{2(d+1)}\alpha-\frac{(\alpha  
+d)(d+1)}{2d^2}\right]  - \frac{d+2}{d}\alpha + 1} \\ \quad  \displaystyle{
\quad =  \frac{\alpha (d^3+d^2 - 2d-1) - d(d+1)^2}{2d^2(d+1)}\alpha^{d-1} -  
\frac{d+2}{d}\alpha + 1}\\ \quad \displaystyle{
\quad \geq  \frac{d (d^3+d^2 - 2d-1) - d(d+1)^2}{2d^2(d+1)}\alpha^{d-1} -  
\frac{d+2}{d}\alpha + 1} \ (\text{since}\ \alpha \geq d\geq 3) \\ \quad \displaystyle{
\quad =  \frac{d^3- 4d-2 }{2d(d+1)}\alpha^{d-1} - \frac{d+2}{d}\alpha + 1} =: B
\end{array}$$
If $\alpha =3$, then $d=3$ and $B= 7/8$. For $\alpha \geq 4$, since $\alpha^{d-1}\geq  
4\alpha $, we further get
$$B \geq \frac{2(d^3- 4d-2 )}{d(d+1)}\alpha - \frac{d+2}{d}\alpha + 1= \frac{d[d(2d-1)-11]-6  
}{d(d+1)}\alpha + 1 >1.$$
Thus we always have $B>0$, which proves (\ref{EB3}).

Now we consider the case when an edge of $\Pcal$ is full, i.e. there are exactly $\alpha +1$  
points on it belonging to $\Acal$. If $\alpha \geq 6$, then $\alpha \geq \frac{4d}{d-1}$ and  
the second condition is satisfied, so we are done. Since $\alpha \geq d$, the left cases are  
$d=4, \ \alpha \leq 5$ and $d=3,\ \alpha =4,5$. In these cases, by Lemma \ref{A2}, $r(S)  
\leq \alpha^{d-2}$. 

If $d=4,\ \alpha \leq 5$, then $\deg K[S] - c+1 \geq \alpha^3 - {\alpha +3 \choose 3} + 6 >  
2\alpha^2\geq 2r(S)$, and by Theorem \ref{A1} we are done.

If $d=3,\ \alpha \leq 5$, let $\tilde{c} = \sharp(M_{\alpha,3}\setminus \Acal)$. Then   
$r(S)\leq \alpha $, and the inequality
$$ \deg K[S] - c+1 = \alpha^2 - {\alpha +2 \choose 2} + \tilde{c} + 4 \geq 2\alpha -1$$
does not hold only in the following situations: $\alpha =3,4,\ \tilde{c}=1,2$ and $\alpha =5, \  
\tilde{c}=1$. By Theorem \ref{A1}, we can restrict ourselves to these situations. By Corollary  
\ref{C3}, we may assume that one deleting point is $(\alpha -1, 1,0)$. Thus, in each case there  
are only few configurations to consider. Using computer, we can check that $r(S)= 2$ if  
$\alpha =3,5$, and $r(S) \leq 3$ if $\alpha =4$. But then $\deg K[S] - c+1 = \alpha^2 - {\alpha  
+2 \choose 2} + \tilde{c} + 4 \geq 2r(S)-1$. Again by Theorem \ref{A1} we are done.
 \hfill $\square$
\end{pf}

{\bf Acknowledgment}.  We would like to thank the referee for valuable comments. This paper was initiated during the visit of the  second  author at Max-Planck Institute for Mathematics in the Sciences (Germany). He would like to thank the MIS for the financial support and the hospitality. All computations  in this paper were done by using the package  
CoCoa \cite{CCA}.

\end{document}